\documentclass[11pt]{article}

\usepackage{amssymb, amsmath, latexsym, mathrsfs, verbatim, calc}

\usepackage[latin1]{inputenc}   

\newtheorem{theorem}{Theorem}

\newtheorem{proposition}[theorem]{Proposition}
\newtheorem{remark}[theorem]{Remark}

\newcommand{\FR}{{\cal F}}

\newcommand{\HR}{{\cal H}}

\newcommand{\MR}{{\cal M}}

\newcommand{\SR}{{\cal S}}

\newcommand{\UR}{{\cal U}}

\newcommand{\WR}{{\cal W}}

\newcommand{\F}{\mathbb{F}}

\newcommand{\N}{\mathbb{N}}

\newcommand{\R}{\mathbb{R}}

\begin{document}

\vskip -0.5cm

\title{Existence of an Optimal Control for Stochastic Control Systems with Nonlinear Cost Functional}
\author{R. Buckdahn$^*$, B. Labed$^\dagger$,  C. Rainer{\footnote{Universit\'e de Bretagne Occidentale, Laboratoire de  Math\'ematiques, CNRS-UMR 6205; 6, av. V. Le Gorgeu, CS 93837, 29238 BREST Cedex 3, France}}, L. Tamer{\footnote{Universit\'e Mohamed Khider-Biskra,
Facult\'{e} des Sciences, D\'{e}partement de Math\'{e}matiques; B.P. 145 Biskra 07000
Alg\'{e}rie}}}
\date{08/12/02}

\maketitle

\vskip -0.4cm

\noindent\textbf{Abstract} We consider a stochastic control problem which is composed of a controlled stochastic differential equation, and whose associated cost functional is defined through a controlled backward stochastic differential equation. Under appropriate convexity assumptions on the coefficients of the forward and the backward equations we prove the existence of an optimal control on a suitable reference stochastic system.
The proof is based on an approximation of the stochastic control problem by a sequence of control problems with smooth coefficients, admitting an optimal feedback control. The quadruplet formed by this optimal feedback control and the associated solution of the forward and the backward equations is shown to converge in law, at least along a subsequence. The convexity assumptions on the coefficients then allow to construct from this limit an admissible control process which, on an appropriate reference stochastic system, is optimal for our stochastic control problem.



\section{Introduction}
In this paper we consider a controlled decoupled forward-backward stochastic differential system of the type
\begin{equation}
\label{equa}
\left\{\begin{array}{l}
dX^{u}_s=b(X^{u}_s,u_s)ds+\sigma(X^{u}_s,u_s)dW_s,\\
dY^{u}_s=-f(X^{u}_s,Y^{u}_s,Z^{u}_s,u_s)ds+Z^{u}_sdW_s+dM^{u}_s, \; s\in [t,T],\\
\langle M^{u},W\rangle_s=0, \; s\in [t,T],\\
X^{u}_t=x,\; Y^{u}_T=\Phi(X^{u}_T), M^{u}_t=0,\\
\end{array}
\right.\end{equation}
where, on a given filtered probability space $(\Omega,\FR,P,\F)$, $W$ is a $d$-dimensional Brownian motion with respect to the not necessarily Brownian filtration $\F$, $X^{u},Y^{u},Z^{u}$ are square integrable adapted processes and $M^{u}$ a square integrable martingale that is orthogonal to $W$. The control problem consists in minimizing the cost functional
$Y^{u}_t$ over all adapted control processes $u$ taking their values in a fixed compact metric space $U$ :
\begin{equation}
\label{original}
V(t,x)=\mbox{essinf}_{u}Y^{u}_t.
\end{equation}
If the driver $f$ of the backward stochastic differential equation (BSDE) doesn't depend on $(y,z)$ the cost functional takes the particular form
\begin{equation}
\label{catherine}
Y^{u}_t=E\left[\Phi(X^{u}_T)+\int_t^Tf(X^{u}_s,u_s)ds~\vert ~{\FR}_t\right],
\end{equation}
reducing the above control problem to the classical one, which has been well studied by a lot of authors; the reader is referred, for instance, to the book \cite{fs} of Fleming and Soner and the references therein. This classical stochastic control problem and its relation with Hamilton-Jacobi-Bellman equations has been generalized by Peng in \cite{peng92}: He characterizes the value function $V(t,x)$ of the stochastic control problem (\ref{original}) as the unique viscosity solution of the associated Hamilton-Jacobi-Bellman equation (the reader is also referred to \cite{peng} for an approach based on BSDE methods). Moreover, in \cite{peng93} he derives a necessary condition for the optimality of a stochastic control for (\ref{original}), given in form of a maximum principle. At the same period,  motivated by applications in econometrics and mathematical finance, Duffie and Epstein have introduced in \cite{de} a similar cost function  called ``recursive utility''; their cost 
 functional $Y^{u}_t$ corresponds to the solution of the above BSDE in (1) if the driver $f$ is supposed not to depend on $z$. For further contributions concerning this stochastic control problem and its applications the reader is referred to \cite{epq}, \cite{dz}, \cite{peng}, \cite{bi} and the references cited therein.

The objective of our present work is to investigate the question of the existence of an optimal control for the problem (\ref{original}).
In the case of a usual stochastic control problem with classical cost functionals (\ref{catherine}) it is well known that, in general, the existence of an optimal control can be got only in the larger class of relaxed controls. However, by El Karoui, Nguyen and Jeanblanc \cite{enj} and by Haussmann and Lepeltier \cite{hl} it has been proved that an optimal control in the original ``strong'' sense exists under some convexity assumptions. We prove an analogous result for the more general case of the controlled forward-backward system of type (\ref{equa}). For this, we approximate (\ref{equa}) by a sequence of stochastic control systems ($\delta$) with smooth coefficients. Since the Hamilton-Jacobi-Bellman equation associated with a stochastic control system with smooth coefficients and strictly elliptic diffusion coefficient admits a smooth solution, it is possible to determine explicitly an optimal feedback control $u^\delta$ by applying a verification theorem. We prove that 
the value functions $V^\delta$ of these approximating systems converge to the value function $V$ of the original problem (\ref{original}). Further, the solution $(X^\delta,Y^\delta,Z^\delta)$ of the approximating control system with the optimal feedback control $u^\delta$ can be approached by a sequence of simple forward equations satisfied by some couple $(\bar X^\delta,\bar Y^\delta)$ and depending on a couple of bounded controls $(\bar Z^\delta,\bar u^\delta)$. This point of view enables us to apply the well known theory
of stochastic controlled forward-systems: We claim that the sequence 
$(X^\delta,Y^\delta)$  converges in law to some couple $(\bar X,\bar Y)$, at least along a subsequence,  and that  this couple satisfies a stochastic differential equation depending on a relaxed control, that is optimal for (\ref{original}). Finally, under a suitable convexity assumption generalizing that of El Karoui, Nguyen and Jeanblanc \cite{enj}, we deduce that this optimum is also attaint by a control in the ``strong sense'', i.e., an admissible control process $u$ defined on an appropriate reference stochastic system $(\Omega,\FR,P,\F,W)$. Moreover, we discuss  our convexity condition and compare it with that of \cite{enj}.

Our paper is organized as follows:
In Section 2, we state the problem and give the main result. Moreover we introduce necessary notations and recall known results which will be used in what follows. Section 3 is devoted to the study of the approximating control problem and the associated Hamilton-Jacobi-Bellman equation: we prove that for the approximating control problem ($\delta$) the essential infimum of the cost functionals is attained by a feedback control $u^\delta$, and that the value function $V^\delta(t,x)$ converges to $V(t,x)$. Finally, in Section 4, we prove our main result concerning the convergence in law of the couples $(X^\delta,Y^\delta)$ along some subsequence $\delta\searrow 0$, and the existence of an optimal control of  (\ref{original}) under an appropriate convexity assumption. Finally, our convexity condition is compared with that of \cite{enj}.

\section{Notations, preliminaries and main theorem}
\label{intro}

Let  $T>0$ be a finite time horizon and $U$ a compact metric space. We  call a reference stochastic system $\nu^0$ a complete probability space $(\Omega^0,\FR^0,P^0)$ endowed with a filtration $\F^0$ satisfying the usual assumptions (i.e. $\F^0$ is right-continuous and ${\FR}^0_0$ contains all $P^0$-null sets in $\FR$) and, according to our needs, with one or two independent $d$-dimensional Brownian motions: $\nu^0=(\Omega^0,\FR^0,P^0,\F^0,W^0)$ or $\nu^0=(\Omega^0,\FR^0,P^0,\F^0,W^0,B^0)$. This reference stochastic system will vary all along this work. On $\nu^0$, we now introduce the following spaces of
processes:\\
For all dimension $m\in\N^*$ and any $t\in[0,T]$,
\begin{itemize}
\item
$\SR^2_{\nu^0}(t,T;\R^m)$ will denote the set of $\R^m$-valued, $\F^0$-adapted, continuous processes $(\Psi_s,s\in[t,T])$ that satisfy $E[\sup_{t\leq s\leq T}|\Psi_s|^2]<\infty$,
\item
$\HR^2_{\nu^0}(t,T;\R^m)$ is the set of $\R^m$-valued, $\F^0$-predictable processes $(\Psi_s,s\in[t,T])$ that satisfy
$E[\int_t^T|\Psi_s|^2ds]<\infty$,
\item ${\cal M}^2_{\nu^0}(t,T;\R^m)$ denotes the set of all $\R^m$-valued, square integrable c\`{a}dl\`{a}g martingales $M=(M_s)_{s\in[t,T]}$ with respect to $\F^0$, with $M_t=0$,
\item $\UR_{\nu^0}(t)$ denotes the set of admissible controls, i.e.   the set of $\F^0$-progressively measurable processes $(u_s,s\in[t,T])$ with values in $U$.
\end{itemize}

\noindent Let us now fix some initial reference stochastic system $\nu=(\Omega,\FR,P,\F,W)$. For any initial condition $(t,x)\in[0,T]\times\R^d$ and any admissible control $u:=(u_s,s\in[t,T])\in\UR_\nu(t)$, we consider the following decoupled forward-backward stochastic differential system:
\begin{equation}
\label{s}
\left\{ \begin{array}{l}
dX_s=b(X_s,u_s)ds+\sigma(X_s,u_s)dW_s,\\
dY_s=-f(X_s,Y_s,Z_s,u_s)ds+Z_sdW_s+dM_s, \; s\in [t,T],\\
X_t=x,\; Y_T=\Phi(X_T), M_t=0, \langle M,W\rangle=0,\\
(X,Y,Z,M)\in\SR_{\nu}(t,T;\R^d)\times\SR_\nu(t,T;\R)\times\HR_\nu(t,T;\R^d)\times\MR^2_{\nu}(t,T;\R^d),
\end{array}
\right.
\end{equation}
where\\
\centerline{$b:\R^d\times U\rightarrow\R^d$, $\sigma:\R^d\times U\rightarrow\R^{d\times d}$,
$f:\R^d\times\R\times\R^d\times U\rightarrow\R$ and $\Phi:\R^d\rightarrow\R$}\\
\noindent satisfy the following assumptions (see, e.g., \cite{bi} or \cite{bl}):
\begin{enumerate}
\item
\begin{itemize}
\item $b$ and $\sigma$ are bounded by some constant $M>0$,
\item for all $x\in\R^d$, $b(x,\cdot)$ and $\sigma(x,\cdot)$ are continuous in $v\in U$,
\item there exists some $C>0$, such that, for all $x,x'\in\R^d$ and $v\in U$,
    \[
    |b(x,v)-b(x',v)|+|\sigma(x,v)-\sigma(x',v)|\leq C|x-x'|.\]
    \end{itemize}
\item
\begin{itemize}
\item $f$ and $\Phi$ are bounded,
\item for all $(x,y,z)\in\R^d\times\R\times\R^d$, $f(x,y,z,\cdot)$ is continuous in $v\in U$,
\item for all $x,x'\in\R^d$, $y,y'\in\R$, $z,z'\in\R^d$ and $v\in U$,
\[
|\Phi(x)-\Phi(x')|+|f(x,y,z,v)-f(x',y',z',v)|\leq C(|x-x'|+|y-y'|+|z-z'|).
\]
\end{itemize}
\end{enumerate}
Under the assumptions 1.-2., the system (\ref{s}) has an unique  solution
 $(X^{t,x,u},Y^{t,x,u},Z^{t,x,u},M^{t,x,u})\in\SR_{\nu}(t,T;\R^d)\times\SR_\nu(t,T;\R)\times\HR_\nu(t,T;\R^d)\times\MR^2_{\nu}(t,T;\R^d)$, and
 \[ J(t,x,u):=Y^{t,x,u}_t\]
 is well defined for all $(t,x)\in[0,T]\times\R^d$ and $u\in\UR_\nu(t)$ (see, for instance, \cite{bi}). We set
 \[ V(t,x)=\mbox{essinf}_{u\in\UR_\nu(t)}J(t,x,u).\]
 Further it is known (see, e.g., \cite{bi}, \cite{bl}) that the a priorily random field $V(t,x)$ possesses a continuous, deterministic version (with which we identify it) and solves in viscosity sense the following Hamilton-Jacobi-Bellman equation:
 \begin{equation}
 \label{hjb}
 \left\{
 \begin{array}{l}
 \displaystyle\frac{\partial}{\partial t}V(t,x)+\inf_{v\in U}H(x,V(t,x),DV(t,x),D^2V(t,x),v)=0,\; (t,x)\in[0,T]\times\R^d,\\
 V(T,x)=\Phi(x),\; x\in\R^d,
 \end{array}\right.
 \end{equation}
 with, for all $(x,y,p,A,v)\in\R^d\times\R\times\R^d\times\mathbb{S}^{d}\times U$,
 \[ H(x,y,p,A,v)=\left\{
 \frac 12\mbox{tr}\left((\sigma\sigma^*)(x,v)A\right)+b(x,v)p
+f(x,y,p\sigma(x,v),v)\right\}\; ,\]
where $\mathbb{S}^{d}$ is the space of the symmetric matrices in $\R^{d\times d},$ and $DV$ and $D^2V$ represent, respectively, the gradient and the Hessian matrix of $V$.\\

\noindent 
Further we consider the following assumption:
\[(H)\left\{ \begin{array}{l}
 \mbox{For all $(x,y)\in\R^d\times\R$, there exists a compact set $A$ in $\R^d\times\R^d\times U$, with}\\
\;\;\;\; A\supset\left\{\;(\sigma^*(x,v)w,0,v)|v\in U,w\in \R^d \mbox{ s.t. } |\sigma^*(x,v)w|\leq K\;\right\}\\
 \mbox{and such that the following set is convex:}\\
\;\;\;\;\{ ((\Sigma\Sigma^*)(x,y,z,\theta,v),\beta(x,y,z,\theta,v))|(z,\theta,v)\in A\}\; ,
\end{array}\right.\]
where, for all $(x,y,z,\theta,v)\in \R^d\times\R\times\R^d\times \R^d\times U$, we have set
\[ \Sigma(x,y,z,\theta,v)=\left(\begin{array}{ll}
\sigma(x,v)&0\\
z^*&\theta^*
\end{array}\right)\;\; \mbox{ and }\;\;
\beta(x,y,z,\theta,v)=\left(\begin{array}{c}
b(x,v)\\
-f(x,y,z,v)
\end{array}\right)\; .\]
The main result of this paper can now be stated:

\begin{theorem}
Under assumption (H), for all $(t,x)\in[0,T]\times\R^d$, there existe a reference stochastic system $\bar\nu$ and an admissible control $\bar u\in\UR_{\bar\nu}(t)$ that is optimal for (\ref{original}), i.e. if 
$(\bar X^{t,x,\bar u},\bar Y^{t,x,\bar u},\bar Z^{t,x,\bar u},\bar M^{t,x,\bar u})
\in\SR_{\bar\nu}(t,T;\R^d)\times\SR_{\bar\nu}(t,T;\R)\times\HR_{\bar\nu}(t,T;\R^d)\times\MR^2_{\bar\nu}(t,T;\R^d)$ is the solution of (\ref{equa}) on $\bar\nu$, we have
\[ \bar Y^{t,x,\bar u}_t=V(t,x)=\mbox{essinf}_{u\in\UR_{\bar\nu}}(t,x,u).\]
\end{theorem}

\noindent We prove this theorem in chapter \ref{convergence}. It is included in Theorem \ref{theo}. Chapter \ref{approxi} prepares this proof by the introduction of an approximating control problem.

\section{An Approximating Control problem}
\label{approxi}

The purpose of this section is to study a sequence of stochastic control problems for which we can explicitly determine an optimal feedback control process and whose value functions converge to that of our original problem. For this end we have to approximate the coefficients of our original control problem by smooth coefficients.

For an arbitrary dimension $m\geq 1$ we let $\varphi :\mathbb{R}^{m}\mathbb{\rightarrow R}$ be a non-negative smooth function on the Euclidean space $\mathbb{R}^{m}$ such that its support is included in the unit ball of $\mathbb{R}^{m}$ and $\int_{%
\mathbb{R}^{m}}\varphi \left( \xi \right) d\xi =1.$ For Lipschitz functions $l:\mathbb{R}^{m}\rightarrow \mathbb{R}$ we set
\begin{equation*}
l_{\delta }\left( \xi \right) =\delta ^{-m}\int_{\mathbb{R}^{m}}l\left( \xi
-\xi ^{^{\prime }}\right) \varphi \left( \delta ^{-1}\xi ^{^{\prime
}}\right) d\xi ^{^{\prime }},\quad \xi\in\mathbb{R}^{m},\, \delta >0.
\end{equation*}%
Then we can easily show that
\begin{equation}
\left\vert l_{\delta }\left( \xi \right) -l\left( \xi \right) \right\vert
\leq C_l\delta,\, \, \, \left\vert l_{\delta }\left( \xi \right) -l_{\delta' }\left( \xi \right) \right\vert
\leq C_l\vert\delta-\delta'\vert,\, \,  \mbox{ for all }\xi \in \mathbb{R}^{m},\, \delta,\delta' >0,  \label{EX1}
\end{equation}
where $C_{l}$ denotes the Lipschitz constant of $l.$

For each $\delta\in(0,1]$ we denote by $b_{\delta },\sigma _{\delta
},f_{\delta }$ and $\Phi _{\delta }$ the mollifications of the functions $b,\sigma ,f$ and $%
\Phi,$ respectively, introduced in Section \ref{intro}, with $l=b\left(
.,v\right) ,$ $\sigma \left( .,v\right) ,$ $f\left( .,v\right) $ and $\Phi
\left( .\right) .$  We emphasize that the estimate (\ref{EX1}) with $l=b\left( .,v\right) ,$ $%
\sigma \left( .,v\right) ,$ $f\left( .,v\right) $ does not depend on $v\in U.$

\medskip

Let us  now fix an arbitrary $\delta\in(0,1]$ and consider the following Hamilton-Jacobi-Bellman equation
\begin{equation}
\label{hjbdelta}
\left\{
\begin{array}{l}
\displaystyle\frac{\partial}{\partial t}V^{\delta}\left( t,x\right)+\underset{v\in U}{
\inf }H^{\delta }\left(x,(V^{\delta },DV^{\delta },D^{2}V^{\delta
})(t,x),v\right) =0,\ \left( t,x\right) \in \left[ 0,T\right] \times \mathbb{R}
^{d},\\
V^{\delta}\left( T,x\right) =\Phi _{\delta }(x),\ x\in \mathbb{R}^{d},
\end{array}%
\right.
\end{equation}%
with the Hamiltonian
\begin{equation*}
H^{\delta }\left( x,y,p,A,v\right)=\frac{1}{2}\left(
\mbox{tr}\left((\sigma _{\delta }\sigma _{\delta }^{\ast })\left( x,v\right)+\delta^2I_{\R^d} \right)A\right)+
b_{\delta}\left( x,v\right)p+f_{\delta }\left(x,y ,p\sigma_{\delta }\left( x,v\right) ,v\right),
\end{equation*}
for $\left(x,y,p,A,v\right) \in \mathbb{R}^{d}\times \mathbb{R}\times\mathbb{R}^d\times\mathbb{S}^d\times U.$ Since the Hamiltonian is smooth and $(\sigma _{\delta }\sigma _{\delta }^{\ast })\left( x,v\right)+\delta^2I_{\R^d}$ is strictly elliptic, we can conclude that the unique bounded continuous viscosity solution $V^\delta$ of the above equation belongs to $C^{1,2}_b([0,T]\times \R^d)$. For this we can apply the regularity results by Krylov \cite{Krylov1} (see the Theorems 6.4.3 and 6.4.4 in \cite{Krylov1}). This regularity properties of the solution $V^\delta$ and the compactness of the control state space $U$ allow to find a measurable function $v^\delta:[0,T]\times\mathbb{R}^d\rightarrow U$ such that, for all $(t,x)\in[0,T]\times\mathbb{R}^d$,
$$H^{\delta }\left(x,(V^{\delta },DV^{\delta },D^{2}V^{\delta
})(t,x),v^\delta(t,x)\right)=\underset{v\in U}{
\inf }H^{\delta }\left(x,(V^{\delta },DV^{\delta },D^{2}V^{\delta
})(t,x),v\right).$$
Let us fix now an arbitrary initial datum $(t,x)\in[0,T]\times\mathbb{R}^d$ and consider the stochastic equation
\begin{equation}
\label{sdedelta}
\left\{
\begin{array}{l}
dX_{s}^{\delta}=b_{\delta }\left( X_{s}^{\delta},v^{\delta
}\left( s,X_{s}^{\delta }\right) \right) ds+\sigma _{\delta }\left(
X_{s}^{\delta},v^{\delta}\left( s,X_{s}^{\delta }\right) \right)
dW_{s}+\delta dB_{s},\, s\in[t,T],\\
X_{t}^{\delta}=x.
\end{array}
\right.
\end{equation}
Since the coefficients $b_{\delta }\left( x',v^{\delta}\left( s,x'\right)\right)$ and $\sigma _{\delta }\left(x',v^{\delta}\left( s,x'\right) \right)$ are measurable and bounded in $(s,x')$ and the matrix $(\sigma _{\delta }\sigma _{\delta }^{\ast })\left( x',v^\delta(s,x')\right)+\delta^2I_{R^d}$ is strictly elliptic, uniformly with respect to $(s,x')\in[t,T]\times \mathbb{R}^d$, we get from Theorem 1 of Section 2.6 in \cite{Krylov2} the existence of a weak solution, i.e., there exists some reference stochastic system
 $(\Omega^{\delta},{\cal F}^\delta,P^\delta,\F^\delta,W^\delta,B^\delta)$  and an $\mathbb{F}^\delta$-adapted continuous process $X^\delta=(X^\delta_s)_{s\in[t,T]}$ such that, $P^\delta$-a.s.,
\begin{equation*}
\left\{
\begin{array}{l}
dX_{s}^{\delta}=b_{\delta }\left( X_{s}^{\delta},v^{\delta
}\left( s,X_{s}^{\delta }\right) \right) ds+\sigma _{\delta }\left(
X_{s}^{\delta},v^{\delta}\left( s,X_{s}^{\delta }\right) \right)
dW^\delta_{s}+\delta dB^\delta_{s},\, s\in[t,T],\\
X_{t}^{\delta}=x.
\end{array}
\right.
\end{equation*}
For an arbitrarily given admissible control $u\in\UR_{\nu^\delta}(t)$, let $X^{\delta,u}$ denote the unique $\mathbb{F}^\delta$-adapted continuous solution of the equation
\begin{equation*}
\left\{
\begin{array}{l}
dX_{s}^{\delta,u}=b_{\delta }\left( X_{s}^{\delta,u},u_s\right) dt+\sigma _{\delta }\left(
X_{s}^{\delta,u},u_s\right)dW^\delta_{s}+\delta dB^\delta_{s},\, s\in[t,T],\\
X_{t}^{\delta,u}=x.
\end{array}
\right.
\end{equation*}
We associate the backward equation
\begin{equation}
\label{bsdedelta}
\left\{
\begin{array}{l}
dY_s^{\delta,u}=-f_{\delta}(X_{s}^{\delta,u},Y_{s}^{\delta,u},Z_{s}^{\delta,u},u_s)ds+Z_{s}^{\delta,u}dW_s^\delta+U_s^{\delta,u} dB_s^\delta+dM_s^{\delta,u},\ s\in[t,T ],\\
Y_T^{\delta,u}=\Phi_\delta(X_{T}^{\delta,u}),\\
(Y^{\delta,u},Z^{\delta,u},U^{\delta,u})\in{\cal S}^2_{\nu^\delta}(t,T;\mathbb{R}) \times {\cal H}^2_{\nu^\delta}(t,T;\mathbb{R}^d)\times {\cal H}^2_{\nu^\delta}(t,T;\mathbb{R}^d), \\
M^{\delta,u}\in{\cal M}^2_{\nu^{\delta}}(t,T;\mathbb{R}^d) \mbox{ is orthogonal to } W^\delta \mbox{ and to } B^\delta.
\end{array}
\right.
\end{equation}
In analogy to our original stochastic control problem we define the cost functionals for our approximating control problem with the help of the solution of (\ref{bsdedelta}):
$$J^\delta(u):=Y^{\delta,u}_t,\, u\in\UR_{\nu^\delta}(t).$$
We shall now identify the solution $V^\delta$ of the Hamilton-Jacobi-Bellman equation (\ref{hjbdelta}) as the value function of our approximating control problem and give an estimate of the distance between the value function $V^\delta$ and that of our original control problem:

\begin{proposition}
\label{prop}
1) Under our standard assumptions we have
$$J^{\delta}(u^\delta)=V^\delta(t,x)=\mbox{essinf}_{u\in\UR_{\nu^\delta}(t)}J^{\delta}(u),$$
where $u^\delta_s:=v^\delta(s,X^\delta_s),\, s\in[0,T],$ is an admissible control from $\UR_{\nu^\delta}(t)$.\\
2) Again under our standard assumptions, there is some constant $C$ only depending on the Lipschitz constants of the functions $\sigma,b,f$ and $\Phi$ such that,
$$\vert V^\delta(t,x)-V(t,x)\vert\le C\delta^{1/2},\, \mbox{ for all } (t,x)\in[0,T]\times \mathbb{R}^d \mbox{ and for all } \delta>0.$$
Moreover, again for some constant $C$ which only depends on the Lipschitz constants and the bounds of the functions $\sigma,b,f$ and $\Phi$ but is independent of $\delta>0$, the following holds true for all $t,t'\in[0,T]$ and $x\in\mathbb{R}^d$:
\begin{equation}
\label{estvd}
\begin{array}{rcl}
\vert V^\delta (t,x)\vert+\vert DV^\delta (t,x)\vert & \le & C,  \\
\vert V^\delta(t,x)-V^\delta(t',x)\vert &\le & C(1+\vert x\vert)\vert t-t'\vert.
\end{array}
\end{equation}

\end{proposition}

{\bf Proof}: 
1) As it is well known that $V^\delta(t,x)=\mbox{essinf}_{u\in\UR_{\nu^\delta}(t)}J^{\delta}(u)$
(see, e.g., \cite{bi}) it only remains to show that $J^{\delta}(u^\delta)=V^\delta(t,x)$. For this end we observe that from the uniqueness of the solution of the controlled forward equation with control process $u^\delta$ it follows that $X^{\delta,u^\delta}=X^{\delta}.$ Moreover, let
$$Y^\delta_s=V^\delta(s,X^{\delta}_s),\,  Z^\delta_s=DV^\delta(s,X^{\delta}_s)\sigma_\delta (X^\delta_s,u^\delta_s), \, U^\delta_s=\delta DV^\delta(s,X^{\delta}_s),\, \,  s\in[t,T],$$
and notice that the triplet $(Y^\delta, Z^\delta,  U^\delta)$ belongs to ${\cal S}^2_{\nu^\delta}(t,T;\mathbb{R}) \times {\cal H}^2_{\nu^\delta}(t,T;\mathbb{R}^d)\times {\cal H}^2_{\nu^\delta}(t,T;\mathbb{R}^d)$.
Then we obtain from It\^{o}'s formula (recall that $V^\delta\in C^{1,2}_b([0,T]\times\mathbb{R}^d)$), combined with the Hamilton-Jacobi-Bellman equation satisfied by $V^\delta$ and the definition of the feedback control $v^\delta(s,x')$, that  $(Y^\delta, Z^\delta,  U^\delta)$ satisfies BSDE (\ref{bsdedelta}) for $u=u^{\delta}$. From the uniqueness of the solution of BSDE(\ref{bsdedelta}) it then follows that $(Y^{\delta,u^\delta},Z^{\delta, u^\delta},U^{\delta,u^\delta})$ $=(Y^{\delta},Z^{\delta},U^{\delta})$ and $M^{\delta,u^\delta}=0.$ Thus, from the definition of $Y^{\delta}$ we get, in particular, that
\[ Y^{\delta,u^\delta}_t=Y^{\delta}_t=V^{\delta}(t,x).\]

\medskip
2) Let $\delta'\in(0,T]$ and $(t',x')\in[0,T]\times \mathbb{R}^d$. Working on the reference stochastic system we have introduced for our arbitrarily fixed $\delta>0$ and $(t,x)\in[0,T]\times \mathbb{R}^d$ at the beginning of this section,
we let $X^{\delta',t',x',u^\delta}\in{\cal S}^2_{\nu^\delta}(t',T;\mathbb{R}^d)$ denote the unique solution of the forward equation
\begin{equation*}
\left\{
\begin{array}{l}
dX_{s}^{\delta',t',x',u^\delta}=b_{\delta'}\left( X_{s}^{\delta',t',x',u^\delta},u^\delta_{s}\right) ds+\sigma _{\delta'}\left(X_{s}^{\delta',t',x',u^\delta},u^\delta_{s}\right)dW_{s}^\delta+\delta' dB_{s}^\delta,\, s\in[t',T],\\
X_{t'}^{\delta',t',x',u^\delta}=x'.
\end{array}
\right.
\end{equation*}
We extend this solution process onto the whole interval $[0,T]$ by setting $X_{s}^{\delta',t',x',u^\delta}=x',$ for $s<t'$. Then, by putting

\medskip

$\displaystyle\widetilde{f}^{\delta',t',x',u^\delta}_s=-\bigg(\frac{\partial}{\partial s}V^{\delta'}(s,X^{\delta',t',x',u^\delta}_s)+ \frac12\mbox{tr}\bigg( \left(\left( \sigma_{\delta'}\sigma_{\delta'}^{\ast}\right) (X^{\delta',t',x',u^\delta}_s,u^\delta_s)+ {\delta'}^2I_{R^d}\right)\times$

$\qquad\displaystyle \times D^2V^\delta(s,X^{\delta',t',x',u^\delta}_s)\bigg)+b_\delta(X^{\delta',t',x', u^\delta}_s,u^\delta_s) DV^{\delta'} (s,X^{\delta',t',x',u^\delta}_s)\bigg), \, s\in [t',T],$

\noindent we define a stochastic process in ${\cal H}^2_{\nu^\delta}(t',T;\mathbb{R})$, and by applying It\^{o}'s formula to $V^{\delta'}(s,X_{s}^{\delta,t',x',u^\delta})$ we show that

\medskip

$Y_s^{\delta',t',x'}=V^{\delta'}(s,X_{s}^{\delta',t',x',u^\delta}),\, Z_s^{\delta',t',x'}=DV^{\delta'}(s,X_{s}^{\delta',t',x',u^\delta}) \sigma_{\delta'}(X_{s}^{\delta',t',x',u^\delta},u^\delta_s),$

$U_s^{\delta',t',x'}=\delta' DV^{\delta'}(s,X_{s}^{\delta',t',x',u^\delta}),\, M_{s}^{\delta', t',x'}=0,\, s\in[t',T],$

\medskip

\noindent is the unique solution of the BSDE
\begin{equation}
\label{deltap}
\left\{
\begin{array}{l}
dY_s^{\delta',t',x'}=-\widetilde{f}^{\delta',t',x',u^\delta}_sds+Z_s^{\delta',t',x'}dW^\delta_s+U_s^{\delta', t',x'}dB^\delta_s +dM_{s}^{\delta',t',x'},\, s\in[t',T],\\
Y_T^{\delta',t',x'}=\Phi_{\delta'}(X_{T}^{\delta',t',x'}),\\
(Y^{\delta',t',x'},Z^{\delta',t',x'},U^{\delta',t',x'})\in{\cal S}^2_{\nu^\delta}(t',T;\mathbb{R}) \times {\cal H}^2_{\nu^\delta}(t',T;\mathbb{R}^d)\times {\cal H}^2_{\nu^\delta}(t',T;\mathbb{R}^d), \\
M^{\delta',t',x'}\in{\cal M}^2_{\nu^{\delta}}(t',T;\mathbb{R}^d) \mbox{ is orthogonal to } W^\delta \mbox{ and to } B^\delta.
\end{array}
\right.
\end{equation}
We want to compare the first component of the solution of (\ref{deltap}) with that of the BSDE
\begin{equation*}
\left\{
\begin{array}{l}
dY_s^{\delta',t',x',u^\delta}=-f_{\delta'}\left(X_{s}^{\delta',t',x',u^\delta},Y_s^{\delta',t',x',u^\delta}, Z_s^{\delta',t',x',u^\delta},u^\delta_s\right)ds\\
\qquad\qquad\qquad +Z_s^{\delta',t',x',u^\delta}dW^\delta_s+U_s^{\delta',t',x',u^\delta} dB^\delta_s +dM_{s}^{\delta',t',x',u^\delta},\, s\in[t',T],\\
\hskip 2mm Y_T^{\delta',t',x',u^\delta}=\Phi_\delta(X_{T}^{\delta,t',x',u^\delta}),\\
(Y^{\delta',t',x',u^\delta},Z^{\delta',t',x',u^\delta},U^{\delta,t',x'},u^\delta)\in{\cal S}^2_{\nu^\delta} (t',T;\mathbb{R}) \times {\cal H}^2_{\nu^\delta}(t',T;\mathbb{R}^d)\times {\cal H}^2_{\nu^\delta} (t',T;\mathbb{R}^d), \\
M^{\delta',t',x',u^\delta}\in {\cal M}^2_{\nu^{\delta}}(t',T;\mathbb{R}^d) \mbox{ is orthogonal to } W^\delta \mbox{ and to } B^\delta.
\end{array}
\right.
\end{equation*}
For this we observe that from the Hamilton-Jacobi-Bellman equation with the classical solution $V^{\delta'}$ it follows that
$$\widetilde{f}^{\delta',t',x',u^\delta}_s\le f_{\delta'}\left(X_{s}^{\delta',t',x',u^\delta},Y_s^{\delta',t',x'}, Z_s^{\delta',t',x'},u^\delta_s\right),\, s\in[t',T].$$
Then the Comparison Theorem for BSDEs yields $Y^{\delta',t',x'}_s\le Y_s^{\delta',t',x',u^\delta},$ $s\in[t',T]$, $P$-a.s., and, consequently, due to 1),
$$V^{\delta'}(t',x')-V^\delta(t,x)=Y^{\delta',t',x'}_{t'}-Y^{\delta,u^\delta}_t\le Y_{t'}^{\delta',t',x',u^\delta} -Y^{\delta,u^\delta}_t,\, P\mbox{-a.s.}$$
On the other hand, SDE and BSDE standard estimates show that there is some generic constant $C$ which depends only on the Lipschitz and the growth constants of the involved functions $b(.,v),\sigma(.,v),f(.,.,.,v)$ and $\Phi(.)$ but  neither on $\delta,\delta'\in(0,1]$ nor on $(t',x'), (t,x)$, such that, with the usual convention $Y_s^{\delta',t',x',u^\delta}=Y_{t'}^{\delta',t',x',u^\delta},$ $Z_{s}^{\delta',t',x',u^\delta}=0,$ $U_{s}^{\delta',t',x',u^\delta}=0,$ and $M_{s}^{\delta',t',x',u^\delta}=0,$ for all $s<t',$

\medskip

$\displaystyle E\bigg[\sup_{s\in[0,T]}\vert Y_s^{\delta',t',x',u^\delta}-Y_s^{\delta,u^\delta}\vert^2$

$\qquad\displaystyle +\int_0^T\left(\vert  Z_s^{\delta',t',x',u^\delta}-Z_s^{\delta,u^\delta}\vert^2 +\vert  U_s^{\delta',t',x',u^\delta}-U_s^{\delta,u^\delta}\vert^2\right)ds+\langle   M^{\delta,t',x',u^\delta} \rangle_T \big\vert{\cal F}^\delta_t\bigg]$

$\le C\left((\delta-{\delta'})^2
+E\big[\sup_{s\in[0,T]}\vert X_s^{\delta',t',x',u^\delta}-X_s^{\delta,u^\delta}\vert^2\vert{\cal F}^\delta_t\big]\right)$

$\le C(\delta-{\delta'})^2+C(1+\vert x\vert^2+\vert x'\vert^2)\vert t-t'\vert+ C\vert x-x'\vert^2$

\medskip

\noindent (Recall that $M^{\delta,v^\delta}=0$). Consequently,

\medskip

$V^{\delta'}(t',x')-V^\delta(t,x)=Y_{t'}^{\delta,t',x'}-Y^{\delta,u^\delta}_t\le Y_{t'}^{\delta',t', x',u^\delta}-Y^{\delta,u^\delta}_t$

$\le C\vert\delta -\delta'\vert+C(1+\vert x\vert+\vert x'\vert)\vert t'-t\vert^{1/2}+ C\vert x-x'\vert,$

\medskip

\noindent and from the symmetry of the argument we get
$$\vert V^{\delta'}(t',x')-V^\delta(t,x)\vert\le C\vert\delta -\delta'\vert+C(1+\vert x\vert+\vert x'\vert)\vert t'-t\vert^{1/2}+ C\vert x-x'\vert.$$
It follows that, in particular,

\medskip

$\vert V^{\delta}(t,x')-V^\delta(t,x)\vert\le C\vert x-x'\vert,$

$\vert V^{\delta}(t',x)-V^\delta(t,x)\vert\le C(1+\vert x\vert)\vert t'-t\vert^{1/2},$

\medskip

\noindent and
$$\vert V^{\delta'}(t,x)-V^\delta(t,x)\vert\le C\vert\delta -\delta'\vert.$$
Moreover, a standard estimate for the BSDE satisfied by $(Y^{\delta,u^\delta},Z^{\delta,u^\delta},U^{\delta,u^\delta},M^{\delta,u^\delta})$ yields the boundedness of $V^{\delta}$, uniformly with respect to $\delta>0$. \\
Therefore, the function $V^\delta$ converges uniformly towards a function $\widetilde{V}\in C_b([0,T]\times \mathbb{R}^d)$, as $\delta$ tends to zero. Since the Hamiltonian $H^\delta$ converges uniformly on compacts to the Hamiltonian of the equation for $V$ it follows from the stability principle for viscosity solutions that $\widetilde{V}$ is a viscosity solution of the same equation as $V$. Thus, from the uniqueness of the viscosity solution within the class of continuous function with at most polynomial growth we get that $\widetilde{V}$ and $V$ coincide. Consequently, $V^{\delta'}$ converges uniformly to $V$, as $\delta'\rightarrow 0$, and from the above estimate of the distance between $V^{\delta'}$ and $V^\delta$ it then follows that

\smallskip

\qquad$\vert V^{\delta}(t,x)-V(t,x)\vert\le C\delta,\mbox{ for all }\delta\in(0,1]\mbox{ and } (t,x)\in[0,T]\times \mathbb{R}^d.$ $\bullet$

\section{Convergence of the Approximating Control Problems}
\label{convergence}
After we have shown that the value function of the approximating problem converges to the value function of the initial problem, we will prove in this  section that there exists a sequence of approximating stochastic controlled systems $(X^{\delta_n},Y^{\delta_n},Z^{\delta_n},u^{\delta_n})$  that converges in law to some controlled system, provided that the couple $(Z^{\delta_n},u^{\delta_n})$ is interpreted as a relaxed control. Then, under some additional convexity condition, we shall find on a suitable reference stochastic system some admissible control - now in the strong sense-, that  is optimal for the initial problem (\ref{original}).\\
We begin this section with the introduction of this additional assumption, then we prove the main result of existence of an optimal control for (\ref{original}) and after we discuss the additional assumption, in particular we compare it with that of \cite{enj}.\\

For all $k>0$, let us denote by $\bar B_k(0)$ the closed ball in $\R^d$ of center 0 and radius $k$. We also introduce the constant K=CM, where
$C>0$ stands here for the constant of the estimates (\ref{estvd}) and $M$  for an upper bound of
$\{ |\sigma(x,v)|, (x,v)\in\R^d\times U\}$.\\
We recall the definition of $(\Sigma,\beta) :\R^d\times\R\times\R^d\times \R^d\times U\rightarrow \R^{(d+1)\times(d+1)}\times
\R^{d+1}$  already introduced in chapter 2:  \\
For all $(x,y,z,\theta,v)\in \R^d\times\R\times\R^d\times \R^d\times U$, we set
\[ \Sigma(x,y,z,\theta,v)=\left(\begin{array}{ll}
\sigma(x,v)&0\\
z^*&\theta^*
\end{array}\right)\;\; \mbox{ and }\;\;
\beta(x,y,z,\theta,v)=\left(\begin{array}{c}
b(x,v)\\
-f(x,y,z,v)
\end{array}\right)\; .\]
We also recall the assumption $(H)$:
\[(H)\left\{ \begin{array}{l}
 \mbox{For all $(x,y)\in\R^d\times\R$, there exists a compact set $A$ in $\R^d\times\R^d\times U$, with}\\
\;\;\;\; A\supset\left\{\;(\sigma^*(x,v)w,0,v)|v\in U,w\in \R^d \mbox{ s.t. } |\sigma^*(x,v)w|\leq K\;\right\}\\
 \mbox{and such that the following set is convex:}\\
\;\;\;\;\{ ((\Sigma\Sigma^*)(x,y,z,\theta,v),\beta(x,y,z,\theta,v))|(z,\theta,v)\in A\}\; .
\end{array}\right.\]

\begin{theorem}
\label{theo}
Suppose that assumption $(H)$ holds and let $(t,x)\in[0,T]\times\R^d$ and $(\delta_n)_{n\in\N}\subset(0,+\infty)$ with $\lim_{n\rightarrow+\infty}\delta_n=0$.
 Then there exists a reference stochastic system \\
 $\bar\nu=(\bar\Omega,\bar\FR,\bar P,\bar\F,\bar W)$, a quadruple
 $(\bar X,\bar Y,\bar Z,\bar M)\in
 \SR^2_{\bar\nu}(t,T;\R^d)\times\SR^2_{\bar\nu}(t,T;\R)\times\SR^2_{\bar\nu}(t,T;\R^d)\times\MR^2_{\bar\nu}(t,T;\R^d)$, with $M$ orthogonal to $W$, and an admissible control $\bar u\in\UR_{\bar\nu}(t)$, such that\\
 1) There is a subsequence of $(X^{\delta_n},Y^{\delta_n})_{n\in\N}$ that converges in distribution to  $(\bar X,\bar Y)$, \\
 2) $(\bar X,\bar Y,\bar Z,\bar M)$ is the solution of the following system:
\begin{equation}
\label{bsde}\left\{
\begin{array}{l}
d\bar X_s=b(\bar X_s,\bar u_s)ds+\sigma( \bar X_s,\bar u_s)d\bar W_s,\\
 d\bar Y_s=-f(\bar X_s,\bar Y_s,\bar Z_s,\bar u_s)ds+\bar Z_sd\bar W_s+d\bar M_s,\;\; s\in[t,T]\\
\bar X_t=x,\; \bar Y_T=\Phi(\bar X_T),
\end{array}\right.\end{equation}
3) For all $(t,x)\in[0,T]\times\R^d$, it holds that
\[\bar Y_t=V(t,x)=\mbox{essinf}_{u\in\UR_{\bar\nu}(t)}J(t,x,u)\; ,\]
i.e. the admissible control $\bar u\in\UR_{\bar\nu}(t)$ is optimal for (\ref{bsde}).
\end{theorem}

\noindent{\bf Proof:}
We shall first introduce an auxiliary sequence of forward-systems, for which the convexity assumption $(H)$ appears as the natural argument to guarantee the existence of a subsequence whose solutions converge in law to a couple $(\bar X,\bar Y)$ associated to a control that is optimal for the original control problem.  Then we will show that the initial sequence $(X^{\delta_n},Y^{\delta_n})_{n\in\N}$ and the auxiliary one have the same limits.\\
1) For all $n\in\mathbb{N}$, let $(X^n,Y^n)$ be the solution  on $\nu^{\delta_n}$ of the following controlled forward system:
\begin{equation}
\label{neq}
\left\{
\begin{array}{l}
dX_{s}^n=b\left( X_{s}^{n},u^{\delta_n}_{s}\right) ds+\sigma \left(X_{s}^{n},u^{\delta_n}_{s}\right)dW_{s}^{\delta_n},\, s\in[t,T],\\
dY_s^{n}=-f\left(X_{s}^{n},Y_s^{n}, w_s^{n}\sigma\left(X_{s}^{n},u^{\delta_n}_{s}\right),u^{\delta_n}_s\right)ds
 +w_s^{n}\sigma\left(X_{s}^{n},u^{\delta_n}_{s}\right)dW^{\delta_n}_s\; s\in[t,T],\\
\hskip 2mm X^n_t=x,\;\; Y^n_t=V(t,x),
\end{array}
\right.
\end{equation}
with $w^n_s=DV^{\delta_n}(s,X^{\delta_n}_s)$.\\
We can rewrite the system (\ref{neq}) as follows:
\begin{equation}
\label{nneq}\left\{
\begin{array}{l}
d\chi^n_s=\beta(\chi^n_s,r^n_s)ds+
\Sigma(\chi^n_s,r^n_s)d\WR^n_s, \; s\in[t,T],\\
\chi^n_t=\left(\begin{array}{c}x\\
V(t,x)\end{array}\right),
\end{array}\right.
\end{equation}
with
\[ \chi^n_s=\left(\begin{array}{c}
X^n_s\\
Y^n_s\end{array}\right),\; r^n_s=(w^n_s\sigma(X^{\delta_n}_s),0,u^{\delta_n}_s)\; \mbox{ and }\; \WR^n=\left(\begin{array}{c}W^{\delta_n}\\B^{\delta_n}\end{array}\right).
\]
Since $w^n_s=DV^{\delta_n}(s,X^{\delta_n}_s)$ and $DV^{\delta}$ is bounded by $C$, uniformly in $\delta$ (Proposition \ref{prop}), we can interpret $(r^n_s,s\in[t,T])$ as a control with values in the compact set $A$ of assumption $(H)$.\\
Now, in order to pass to the limit in $n$, we shall as usual embed the controls $r^n$  in the set of relaxed controls, i.e. consider $r^n$ as random variable with values in the space $V$ of all Borel measures $q$ on $[0,T]\times A$, whose projection $q(\cdot\times A)$ concides with the Lebesgue measure. For this, we identify the control process $r^n$ with the random measure
\[ q^n(\omega,ds,da)=\delta_{r^n_s(\omega)}(da)ds,\; (s,a)\in[0,T]\times A,\omega\in \Omega.\]
  From the boundedness of $\{ \left(\Sigma(x,y,z,\theta,v),\beta(x,y,z,\theta,v)\right), (x,y,z,\theta,v)\in\R^d\times\R\times A\}$ and the compactness of $V$ with respect to the topology induced by the weak convergence of measures, we get the tightness of the laws of $(\chi^n,q^n),n\geq 1,$ on $C([0,T];\R^d\times\R)\times V$.
Therefore we can find a probability measure $Q$ on $C([0,T];\R^d\times\R)\times V$ and extract a subsequence -still denoted by $(\chi^n,q^n)$- that converges in law to the canonical process  $(\chi,q)$ on the space $C([0,T];\R^d\times\R)\times V$ endowed with the measure $Q$.\\
 Now, by assumption $(H)$ on the coefficients of system (\ref{nneq}), we can apply the result of  \cite{enj}, that claims that there exists a stochastic reference system $\bar\nu=(\bar\Omega, \bar\FR,\bar P,\bar\F,\bar\WR)$ enlarging $(C([0,T];\R^d\times\R)\times V;Q)$ and an  $\bar\F$-adapted process  $\bar r$ with values in $A$, such that the process $\chi$  is a solution of
 \[\left\{
\begin{array}{l}
d\chi_s=\beta(\chi_s,\bar r_s)ds+
\Sigma(\chi_s,\bar r_s)d\bar\WR_s, \; s\in[t,T],\\
\chi_t=\left(\begin{array}{c}x\\
V(t,x)\end{array}\right),
\end{array}\right.
\]
and has the same law under $\bar P$ as under $Q$. 
Replacing $\Sigma$ and $\beta$ by their definition and setting $\chi=\left(\begin{array}{c}\bar X\\
 \bar Y\end{array}\right)$, $\bar\WR=\left(\begin{array}{c}\bar W\\
 \bar B\end{array}\right)$ and $\bar r=(\bar Z,\bar\theta,\bar u)$, this system is equivalent to
 \[\left\{\begin{array}{l}
 d\bar X_s=b(\bar X_s,\bar u_s)ds+\sigma( \bar X_s,\bar u_s)d\bar W_s,\\
 d\bar Y_s=-f(\bar X_s,\bar Y_s,\bar Z_s,\bar u_s)ds+\bar Z_sd\bar W_s+\bar\theta_sd\bar B_s,\;\; s\in[t,T]\\
\bar X_t=x,\;\bar Y_t=V(t,x).
\end{array}\right.
 \]
 2) It follows from standard estimations that
for some constant $K>0$ and for all $n\in\N$,
\begin{equation}
\label{dn}
\begin{array}{c}
E[\sup_{s\in[t,T]}|X_s^{\delta_n}-X^n_s|^2]\leq K\delta_n,\\
E[\sup_{s\in[t,T]}|Y_s^{\delta_n}-Y^n_s|^2]\leq K\delta_n.
\end{array}
\end{equation}
This implies that if some subsequence of  $(X^n,Y^n)_{n\in\N}$ converges in law, the same holds true for $(X^{\delta_n},Y^{\delta_n})_{n\in\N}$, and the limits have same law.
Further we deduce from (\ref{dn}) and Proposition \ref{prop},  that
$\bar Y_s=V(s,\bar X_s)$ for all $s\in[t,T]$ $\bar P$-a.s.. In particular $\bar Y_T=\Phi(\bar X_T)$ $\bar P$-a.s. .
Thus, if we set $\bar M_s=\int_t^s\bar\theta_rd\bar B_r$, then
$\langle \bar M,\bar W\rangle_s=\int_t^s\bar\theta_rd\langle \bar B,\bar W\rangle_r=0$ and $(\bar Y,\bar Z,\bar M)$ satisfies (\ref{bsde}).\\
3) We have already seen that $\bar Y_s=V(s,\bar X_s)$ for all $s\in[t,T]$ $\bar P$-a.s.  
On the other hand, it is well known that, for the unique bounded viscosity solution $V$ of the Hamilton-Jacobi-Bellman equation (\ref{hjb}),
\[ V(t,x)=\mbox{essinf}_{u\in\UR_{\bar\nu}(t)}J(t,x,u), \; \bar P\mbox{-a.s.}\]
(see e.g. \cite{bi}). Thus assertion 3) of the theorem follows.
 $\bullet$\\

\noindent We state now a proposition that relays assumption $(H)$ to some convexity assumptions concerning the parameters of the initial system (\ref{s}).

\begin{proposition}
\label{h}
1) We suppose that the following assumption holds:
\[(H1)\left\{\begin{array}{l}\mbox{
For all $(x,y)\in\R^d\times\R$, the set}\\
\left\{ \left((\sigma\sigma^*)(x,v),(\sigma\sigma^*)(x,v)w,b(x,v),f(x,y,\sigma^*(x,v)w,v)\right)|v\in U,w\in \R^d
\mbox{ s.t. } |\sigma^*(x,v)w|\leq K\right\} \\
 \mbox{is convex.}\end{array}\right.\]
Then $(H)$ is satisfied.\\
2)
We suppose that $f$ doesn't depend on $z$, i.e., for all $(x,y,z,v)\in\R^d\times\R\times\R^d\times U$, $f(x,y,z,v)=f(x,y,0,v):=f(x,y,v)$.\\
Moreover, we assume that the following assumption is satisfied:
\[ (H2)\left\{\begin{array}{l}
\mbox{For all $(x,y)\in\R^d\times\R$, the set }\\
\{ ((\sigma\sigma^*)(x,v),b(x,v),f(x,y,v))|v\in U\}
\mbox{ is convex}.\hspace{4cm}\end{array}\right.
\]
Then we have $(H)$.\\
\end{proposition}

\begin{remark}
For the case of classical stochastic control problems with $f$ independent of $(y,z)$, we find in $(H2)$ the standard assumption which guaranties the existence of an optimal control on a suitable reference stochastic system (see \cite{enj}).
\end{remark}

\noindent{\bf Proof:} 1) Let us fix $(x,y)\in\R^d\times\R$.
We will show that, under assumption (H1), there exists a set $A\subset\R^d\times\R^d\times U$ such that
\[ \begin{array}{r}
\overline{co}\{ \left((\Sigma\Sigma^*)(x,\sigma^*(x,v)w,0,v),\beta(x,y,\sigma^*(x,v)w,v)\right)|(v,w)\in U\times \R^d \mbox{ s.t. } |\sigma^*(x,v)w|\leq K\}\\
=\{ \left((\Sigma\Sigma^*)(x,y,z,\theta,v),\beta(x,y,z,\theta,v)\right)|(z,\theta,v)\in A\},
\end{array}\]
(where, for any set $E$, $\overline{co}E$ stands for the convex hull of $E$) and that we can choose $A$ compact.\\
For this end, we consider an arbitrarily chosen probability measure $\mu$ on the set $\Gamma=\{ (v,w)\in U\times\R^d|\; |\sigma^*(x,v)w|\leq K\}$.\\
The first step consists in finding a triplet $(\bar z,\bar\theta,\bar v)\in \R^d\times\R^d\times U$
such that
\begin{equation} \label{trip}
\begin{array}{r}
\int_{\Gamma}((\Sigma\Sigma^*)(x,\sigma^*(x,v)w,0,v),\beta(x,y,\sigma^*(x,v)w,v)\mu(dv,dw)\\
=\left((\Sigma\Sigma^*)(x,y,\bar z,\bar\theta,\bar v),\beta(x,y,\bar z,\bar\theta,\bar v)\right).
\end{array}
\end{equation}
Set $\Phi(v,w)=\left((\sigma\sigma^*)(x,v),(\sigma\sigma^*)(x,v)w,b(x,v),f(x,y,\sigma^*(x,v)w,v)\right)$.
As assumption $(H1)$ is supposed to hold true, there exists a couple $(\bar v,\bar w)$ in $\Gamma$ such that
\begin{equation}
\label{vw}
 \int_{\Gamma}\Phi(v,w)\mu(dv,dw)=\Phi(\bar v,\bar w).
 \end{equation}
 We set $\bar z=\sigma^*(x,\bar v)\bar w$.
Now the explicit calculus of $(\Sigma\Sigma^*)(x,y,\sigma^*(x,v)w,0,v)$ shows that, for getting (\ref{trip}), it suffices to find $\bar\theta\in \R^d$ such that
\begin{equation}
\label{theta}
\alpha:=\int_{\Gamma}w^*(\sigma\sigma^*)(x,v)w\mu(dv,dw)-\bar w^*(\sigma\sigma^*)(x,\bar v)\bar w=|\bar\theta|^2.
\end{equation}
But $\alpha$ is non negative, since it can also be written as
\[\alpha =
\int_{\Gamma}\left(\sigma^*(x,v)(w-\bar w)\right)\left(\sigma^*(x,v)(w-\bar w)\right)^*\mu(dv,dw)\geq 0.
\]
Consequently, such $\theta\in\R^d$ satisfying (\ref{theta}) exists.
Further let us rewrite (\ref{theta}) as
\begin{equation}
\label{ztheta}
 \int_{\Gamma}|\sigma^*(x,v)w|^2\mu(dv,dw)=|\bar z|^2+|\bar\theta|^2.
 \end{equation}
Since the support of $\mu$ is included in $\Gamma$, it follows that $\bar z$ and $\bar\theta$ belongs to $\bar B_K(0)$.\\
Now we define $B$ as the set of triplets $(\bar z,\bar\theta,\bar v)$ for which there exist $\bar w\in\R^d$ and a probability measure $\mu$ on $\Gamma$  such that \\
{\it (i)} $\bar z=\sigma^*(x,\bar v)\bar w$, $|\bar z|\leq K$\\
{\it (ii)} the relations  (\ref{vw}) and (\ref{theta}) are satisfied.\\
Consequently, 
\[\begin{array}{r} \overline{co}\left\{ \left( (\Sigma\Sigma^*)(x,\sigma^*(x,v)w,0,v),\beta(x,y,\sigma^*(x,v)w,0,v)\right)|(v,w)\in\Gamma\right\}\\
=\left\{\left( (\Sigma\Sigma^*)(x,y,z,\theta,v),\beta(x,y,z,\theta,v)\right)|(z,\theta,v)\in B\right\}.
\end{array}\]
Le $A$ be the closure of $B$.
 From the boundedneess of $B$ ($\subset\bar B_K(0)\times\bar B_K(0)\times U$) follows the compactness of $A$.
Moreover, since $\left((\Sigma\Sigma^*)(x,y,\cdot,\cdot,\cdot),\beta(x,y,\cdot,\cdot,\cdot)\right)$ is continuous, the convexity of 
$\left\{\left( (\Sigma\Sigma^*)(x,y,z,\theta,v),\beta(x,y,z,\theta,v)\right)|(z,\theta,v)\in B\right\}$ implies that of \\
$\left\{\left( (\Sigma\Sigma^*)(x,y,z,\theta,v),\beta(x,y,z,\theta,v)\right)|(z,\theta,v)\in A\right\}$.
Finally, for all $(\bar v,\bar w)\in U\times\R^d$ with $|\sigma^*(x,\bar v)\bar w|\leq K$, 
$(\sigma^*(x,\bar v)\bar w,0,\bar v)\in B\subset A$ (for $\mu=\delta_{\bar v,\bar w}$).\\

\noindent 2) As for 1), for all probability measure $\mu$ on $ \Gamma$, we can find $\bar v\in U$ such that
\[ \int_{\Gamma}\left((\sigma\sigma^*)(x,v),b(x,v),f(x,y,v)\right)\mu(dv,dw)=
\left((\sigma\sigma^*)(x,\bar v),b(x,\bar v),f(x,y,\bar v)\right).
\]
We will show next that there exists $\bar w\in\R^d$ such that
\begin{equation}
\label{v}
\int_{\Gamma}(\sigma\sigma^*)(x,v)w\mu(dv,dw)
=(\sigma\sigma^*)(x,\bar v)\bar w.
\end{equation}
Indeed, since the matrix $(\sigma\sigma^*)(x,\bar v)$ is semidefinite positive, it can be written as \[(\sigma\sigma^*)(x,\bar v)=T\Lambda T^*,\]
 where
$T\in\R^{d\times d}$ such that $TT^*=T^*T=I_{\R^d}$,
 and
$\Lambda=\left(
\begin{array}{ccc}
\lambda_1&&0\\
&\ldots&\\
0&&\lambda_d
\end{array}\right)\in\R^{d\times d}$, with $\lambda_1\geq\ldots\geq \lambda_d\geq 0$.\\
For $\{e_1,\ldots,e_d\}$ the canonical basis of $\R^d$, we set $f_k=Te_k,\; k\in\{ 1,\ldots,d\}$.
Remark that also $\{ f_1,\ldots,f_d\}$ is an orthonormal basis of $\R^d$.\\
Now let $l\in\{ 0, \ldots,d\}$ be such that
$\lambda_1\geq\ldots\geq \lambda_l>0=\lambda_{l+1}=\ldots=\lambda_d$.
If $l=d$, the matrix $(\sigma\sigma^*)(x,\bar v)$ is invertible and to get (\ref{v}),  we just have to set
\[\bar w=\{(\sigma\sigma^*)(x,\bar v)\}^{-1}\int_{\Gamma}(\sigma\sigma^*)(x,v)w\mu(dv,dw).\]
Else, for all $r\in\{ l+1,\ldots,d\}$, we have
\[ \langle(\sigma\sigma^*)(x,\bar v)f_r,f_r\rangle
=\langle\Lambda e_r,e_r\rangle=0\; ,\]
where $\langle\cdot,\cdot\rangle$ stands for the scalar product in $\R^d$.
This implies that
\[
\begin{array}{rl}
\int_{\Gamma}|\sigma^*(x,v)f_r|^2\mu(dv,dw)=&
\int_{\Gamma}\langle (\sigma\sigma^*)(x,v)f_r,f_r\rangle\mu(dv,dw)\\
=&\langle\left( \int_{\Gamma}(\sigma\sigma^*)(x,v)\mu(dv,dw)\right)f_r,f_r\rangle\\
=&\langle(\sigma\sigma^*)(x,\bar v)f_r,f_r\rangle\\
=& 0.
\end{array}\]
In other words $\sigma^*(x,v)f_r=0$, $\mu$-a.e.,
and
\[ \langle\int_{\Gamma}(\sigma\sigma^*)(x,v)w\mu(dv,dw),f_r\rangle
=\int_{\Gamma}w^*(\sigma\sigma^*)(x,v)f_r\mu(dv,dw)=0\; , l+1\leq r\leq d\; .\]
Hence $\int_{U\times \bar B_C(0)}(\sigma\sigma^*)(x,v)w^*\mu(dv,dw)\in\mbox{span}\{ f_1,\ldots,f_l\}$, and
 the existence of $\bar w\in\R^d$ satisfying (\ref{v}) follows.\\
The end of the proof  proceeds as for statement 1) : we chose some $\bar\theta\in\R^d$ that satisfies (\ref{theta}),
we set $\bar z=\sigma^*(x,\bar v)\bar w$ and prove that $\bar z$ and $\bar\theta$ necessarily belong to $B_K(0)$. Then
we define $A$ as for statement 1) . $\bullet$


\begin{thebibliography}{99}

\bibitem{bi} R. Buckdahn, N. Ichihara, {\it Limit Theoem for Controlled Backward SDEs and Homogenization of Hamilton-Jacobi-Bellman Equations},
Appl. Math. Optim. 51 (2005), pp.1-33.

\bibitem{bl} R. Buckdahn, J. Li, {\it Stochastic differential games and viscosity solutions of Hamilton-Jacobi-Bellman-Isaacs equations},  SIAM J. Control Optim.  47,  no. 1 (2008), pp.444-475.

\bibitem{de} D. Duffie, L. Epstein, {\it Stochastic differential utility}, Econometrica, 60 (1992), pp.353-394.

\bibitem{dz} N. Dokuchaev, X.Y. Zhou, {\it Optimal Controls of Backward Stochastic Differential Equations}, J. Math. Anal. Appl. 238, no. 1 (1999), pp.143-165.

\bibitem{enj} N. El Karoui, D.H. Nguyen, M. Jeanblanc-Piqu\'e, {\it Compactification Methods in the Control of Degenerate Diffusions: Existence of an Optimal Control}, Stochastics, Vol.20 (1987), pp.169-219.

\bibitem{epq} N. El Karoui, S. Peng, M.V. Quenez, {\it Backward stochastic differential equations in finance}, Mathematical Finance, 7  (1997) pp.1-71.

\bibitem{fs} W.~H. Fleming, H.~M.Soner, {\it Controlled Markov processes and viscosity solutions}, 2nd ed., Stochastic Modelling and Applied Probability 25, New York, NY: Springer, 2006.

\bibitem{hl} U.G. Haussmann, J.P. Lepeltier, {\it On the existence of optimal controls}, SIAM J. Control Optim. 28, No.4 (1990), pp.851-902.

\bibitem{Krylov1} N.K.Krylov, {\it Nonlinear Elliptic and Parabolic Equations of Second Order}, Reidel, Dordrecht, 1987.

\bibitem{Krylov2} N.K.Krylov, {\it Controlled Diffusion Processes}, Applications of Mathematics 14, Springer, New York, Heidelberg, Berlin, 1980.

\bibitem{peng92}  S.Peng, {\it A generalized dynamic programming principle and Hamilton-Jacobi-Bellman equation}, Stochastics and
Stochastics Reports, Vol.38 (1992), 119-134.

\bibitem{peng93} S. Peng, {\it Backward Stochastic Differential Equations and
Applications to Optimal Control}, Appl Math Optim 27 (1993), pp.125-144.

\bibitem{peng} S. Peng, {\it BSDE and stochastic optimizations}, Topics in stochastic analysis, Chap.2, S.Yan, S.Peng, S.Fang, and L.Wu, (eds.), Science Press, Beijing, 1997.


\end{thebibliography}
\end{document}